\documentclass[11pt]{article}
\usepackage[T1]{fontenc}
\usepackage{lmodern,microtype}
\usepackage[margin=1.02in]{geometry}
\usepackage{amsmath,amssymb,amsthm,mathtools}
\usepackage{enumitem,xcolor}
\usepackage[colorlinks=true,linkcolor=blue!45!black,citecolor=blue!45!black,urlcolor=blue!45!black]{hyperref}
\hypersetup{pdftitle={The Mezard--Parisi formula for diluted Potts models},pdfauthor={Dmitry Panchenko}}
\numberwithin{equation}{section}
\newtheorem{theorem}{Theorem}[section]
\newtheorem{proposition}[theorem]{Proposition}
\newtheorem{lemma}[theorem]{Lemma}
\newtheorem{corollary}[theorem]{Corollary}
\theoremstyle{definition}
\newtheorem{definition}[theorem]{Definition}
\theoremstyle{remark}
\newtheorem{remark}[theorem]{Remark}
\newcommand{\E}{\mathbb E}
\newcommand{\Prob}{\mathbb P}
\newcommand{\R}{\mathbb R}
\newcommand{\N}{\mathbb N}
\newcommand{\one}{\mathbf 1}
\newcommand{\Poi}{\operatorname{Pois}}
\newcommand{\supp}{\operatorname{supp}}
\newcommand{\Av}{\operatorname{Av}}
\newcommand{\Cov}{\operatorname{Cov}}

\newcommand{\Sp}{\mathsf P}
\newcommand{\Law}{\operatorname{Law}}
\newcommand{\Dir}{\operatorname{Dir}}
\newcommand{\cA}{\mathcal A}

\newcommand{\cG}{\mathcal G}
\newcommand{\cP}{\mathcal P}
\newcommand{\cH}{\mathcal H}

\newcommand{\cc}{\boldsymbol c}
\newcommand{\zz}{\boldsymbol\zeta}
\newcommand{\dd}{\,\mathrm d}
\newcommand{\ip}[2]{\langle #1,#2\rangle}
\newcommand{\nrm}[1]{\lVert #1\rVert}

\newcommand{\MP}{M\'ezard--Parisi}
\title{The \MP\ formula for diluted Potts models}
\author{Dmitry Panchenko\thanks{
Department of Mathematics, University of Toronto,
panchenk@math.toronto.edu.
Partially supported by an NSERC grant.}}
\date{}
\begin{document}
\maketitle
\begin{abstract}
We obtain the Mézard–Parisi formula for the antiferromagnetic Potts model on random uniform hypergraphs, including the Erdős–Rényi graph. The upper bound in the graph case was proved by Contucci, Dommers, Giardinà, and Starr, whose interpolation extends to uniform hypergraphs. Here we prove the matching lower bound for a general class of diluted models with finite spin alphabets. The proof follows closely the corresponding argument for diluted Ising models, so we outline only the adjustments needed for several colors.
\end{abstract}

\section{Introduction}\label{sec:intro}

In this paper, we prove the M\'ezard--Parisi formula for the
antiferromagnetic Potts model at every finite nonnegative
inverse temperature and every nonnegative mean degree, including uniform hyperedges of any fixed arity at least two. It is a direct extension of the \MP\  formula for diluted models with Ising spins in \cite{MPnew}, whose notation and argument we follow. The structure of the proof is essentially unchanged and, as in the Ising case, the principal task is the lower bound.  With several colors, a local magnetization becomes a vector of color probabilities, and products of color indicators replace spin moments. We describe the constructions needed to see what these vectors encode, and explain the adjustments at the points where they enter the proof.

We present this extension separately because the Ising setting makes the reproducibility mechanism especially transparent, with lighter notation. Once that argument is understood, we can concentrate on the adjustments introduced by several colors, which we collect in Section~\ref{sec:proof}.

The preliminary sections state the main results and recall
the Aldous--Hoover representation for several colors.
We introduce the colored multioverlaps and the scalar equality
overlap, then describe the ultrametric geometry determined
by this overlap and explain approximate reproducibility
in this setting. These constructions provide the context
for the proof modifications that follow in Section~\ref{sec:proof}.
In Appendix~\ref{sec:upper}, we extend the interpolation
of \cite{CDGS} to uniform hypergraphs.
Appendix~\ref{app:heps} gives a new short proof of hierarchical
exchangeability \cite{HEPS}, based on the conditional
Dirichlet law of the cluster masses.

Let us briefly mention the main adjustments in the proof compared with the Ising case. The equality overlap is scalar, yet its Gram structure controls fluctuations in every color coordinate; this gives the starting estimate for the induction, simultaneously for all colors. A later adjustment concerns the mixed cavity tests, whose permutations move the assigned colors together with the replica shape. The existing reduction to common-stem
moments resolves this issue when applied before relabeling, as we explain in Section~\ref{sec:graph}. For the simultaneous limiting identities and cavity lower bound, we use the Ising construction in \cite{GGcavInput}, with the adjustments recalled in Section~\ref{sec:joint-input}.

\paragraph{Related work.}

The cavity approach of M\'ezard and Parisi \cite{MP1,MP2}
extended the Parisi ansatz \cite{MPV,Parisi79,Parisi80} to
models with a finite mean number of neighbors per spin.
It predicted both the free energy and the local spin
distributions among states; we refer to \cite{MM,Zreview}
for broader accounts of this approach.

For random graph coloring, Krz\k{a}ka\l a, Pagnani, and
Weigt \cite{KPW} studied the coloring thresholds, the
stability of the one-step replica-symmetry-breaking ansatz,
and the large-$q$ asymptotics. Krz\k{a}ka\l a, Montanari,
Ricci-Tersenghi, Semerjian, and Zdeborov\'a \cite{KMRSZ}
developed the physical picture of clustering and condensation
in random constraint satisfaction problems, including
graph coloring, and predicted Poisson--Dirichlet statistics
for the dominant cluster weights in the condensed phase.
For coloring, this picture was developed in detail by
Zdeborov\'a and Krz\k{a}ka\l a \cite{ZKcoloring}, who
described clustering, condensation, and freezing.
Krz\k{a}ka\l a and Zdeborov\'a \cite{KZPotts} studied the
positive-temperature Potts model and the stability of its
replica-symmetry-breaking solutions. For random uniform
hypergraphs, Gabri\'e, Dani, Semerjian, and Zdeborov\'a
\cite{GDSZ} studied the corresponding $q$-coloring problem
and its phase transitions, including regimes in which
the one-RSB solution is unstable.

For the antiferromagnetic Potts model on the
Erd\H{o}s--R\'enyi graph, Contucci, Dommers, Giardin\`a,
and Starr \cite{CDGS} established an extended variational
principle and the hierarchical interpolation upper bound
used here. Coja-Oghlan and Jaafari \cite{CJ} determined
the free energy at all temperatures in an explicit range
of mean degrees. Coja-Oghlan and Zdeborov\'a \cite{CZhyper}
proved the existence of a condensation transition in
random hypergraph 2-coloring, while Coja-Oghlan,
Krz\k{a}ka\l a, Perkins, and Zdeborov\'a \cite{CKPZ}
gave a variational characterization of the condensation
transition for the Potts antiferromagnet and random graph
coloring. More recently, Biswas, Chen, and Sen \cite{BCS}
proved replica-symmetric formulas for general diluted spin
systems in high-temperature and subcritical regimes,
including Potts interactions.

Several works place these questions in the broader setting
of random factor graphs. Bapst and Coja-Oghlan \cite{BCO}
and Coja-Oghlan and Perkins \cite{CPBP} give sufficient
conditions for the replica-symmetric Bethe free-energy
formula. Coja-Oghlan and Perkins also prove decompositions
into Bethe states \cite{CPStates} and derive a variational
formula for models on random regular graphs \cite{CPBethe}.
For the fully connected Potts spin glass, the Parisi formula
was proved in \cite{PottsMF}, and Bates and Sohn
\cite{BSPotts} reduced its matrix order parameter to a scalar
one under the constraint of balanced color proportions.

The rigorous variational methods have their origins in
Guerra's interpolation bound \cite{Guerra} and the cavity
variational principle of Aizenman, Sims, and Starr
\cite{ASS}. Talagrand \cite{TPF} proved the Parisi formula
for the Sherrington--Kirkpatrick model \cite{SK}, and its
extension to general mixed $p$-spin models was proved in
\cite{PPF}. For diluted models, Franz and Leone \cite{FL}
(see also \cite{PT}) developed the corresponding
interpolation bounds. De Sanctis \cite{DeSanctis} gave
an extended variational principle for the Viana--Bray
model using random multioverlap structures.

The approach we follow uses spin distributions and cavity
equations \cite{Spins}, hierarchical exchangeability
\cite{HEPS}, and the hierarchical representation theorem
\cite{AP}. The earlier reproducibility results for certain
one-RSB models \cite{1RSB} and for finite-RSB models
\cite{rRSB} lead to the proof of the \MP\ formula in the
Ising case in \cite{MPnew}. Another proof of the
M\'ezard--Parisi formula for diluted even-arity Ising
models appeared in \cite{OAI}.\vspace{-4mm}

\paragraph{Acknowledgements.}
The motivation to consider a finite-alphabet version of the
\MP\ formula came from a conversation with Lenka Zdeborov\'a.
In collaboration with ChatGPT Astra, I examined the extension
of the Ising proof and checked the necessary modifications. ChatGPT Astra also found the short argument for conditional independence of the cluster masses and sampled spins that we use in Appendix~\ref{app:heps} to give a new proof of the hierarchical exchangeability in \cite{HEPS}.

\setlength{\parskip}{2pt}

\section{Models and main results}\label{sec:model}

\subsection{Diluted models with a finite alphabet}

We first specify the family of models for which we prove
the lower bound. We use the Poisson clause construction of
\cite{MPnew}, now allowing finitely many colors and several clause
types. We work with the alphabet \([q]=\{1,\ldots,q\}\), where \(q\ge2\), and write
\[
 \Delta_q=\{x\in[0,1]^q:\textstyle\sum_{a=1}^q x_a=1\}
\]
for the simplex of color distributions. Let \(\mathcal I\) be a finite or
countable set of primary clause types. Type \(t\) has arity \(k_t\ge2\), intensity
\(\lambda_t\ge0\), and a random clause
\(\theta_t:[q]^{k_t}\to\R\).
We denote its independent copies by \(\theta_{t,\nu}\).
For each type, we sample a count $M_{t,N}\sim\Poi(\lambda_tN)$.
For each of its clauses, we choose the $k_t$ site indices
independently and uniformly from $[N]$, allowing repetitions.
All counts, indices, and clauses are independent. We then define
\begin{equation}\label{eq:model}
 H_N(\sigma)=\sum_{t\in\mathcal I}\sum_{\nu=1}^{M_{t,N}}
       \theta_{t,\nu}(\sigma_{i_{t,\nu,1}},\ldots,
                     \sigma_{i_{t,\nu,k_t}}),\qquad
 F_N=\frac1N\E\log\sum_{\sigma\in[q]^N}e^{H_N(\sigma)}.
\end{equation}
Writing \(J_t^0=\max_{\boldsymbol a\in[q]^{k_t}}
|\theta_t(\boldsymbol a)|\), we assume the integrability condition
\begin{equation}\label{eq:primary-moments}
 \sum_{t\in\mathcal I}\lambda_t k_t^2
              \bigl(\E e^{sJ_t^0}-1\bigr)<\infty
                      \qquad\text{for every }s>0.
\end{equation}
We use this sufficient condition to make the Hamiltonian
and cavity fields well defined and to carry over the analytic
estimates of the Ising proof. In particular, the exponential
moments control the cavity weights and the Lipschitz bounds needed
when we replace states by their cluster barycenters in the cavity
equations; we use this step in Section~\ref{sec:transfer}, following
\cite[Section~7]{MPnew}. The weight $k_t^2$ accounts for both the
intensity $\lambda_t k_t$ of clauses incident to an added site and
the $k_t-1$ incoming arguments over which these bounds are summed.
The same assumption supplies the moment and truncation bounds
used in the joint limiting construction in
Section~\ref{sec:joint-input} and the variational passage in
Section~\ref{sec:bridge}. It is automatic for finitely many bounded
clause laws, including the Potts interactions considered below.

We may assume that each clause law is invariant under permutations
of its argument positions. This imposes no restriction on the
model: independently and uniformly permuting the argument positions of every
clause is equivalent in distribution to permuting its independent uniform
site indices. This symmetrization concerns the positions
within a clause, not the color labels. No color-permutation
symmetry is assumed. From now on, \(\theta_t\) denotes this
symmetrized clause law, including in the trial functional below.
This convention allows us to place the added site's color in the last
argument of every incident cavity clause.

A cavity clause is evaluated at independent input colors whose
distributions are probability vectors. We therefore use the same
extension convention as in \cite{MPnew}: we average the Boltzmann
factor over these colors. For $x_1,\ldots,x_{k_t}\in\Delta_q$, we set
\begin{equation}\label{eq:extension}
 e^{\theta_t(x_1,\ldots,x_{k_t})}
 =\sum_{a_1,\ldots,a_{k_t}\in[q]}
      e^{\theta_t(a_1,\ldots,a_{k_t})}
                             \prod_{j=1}^{k_t}x_{j,a_j}.
\end{equation}
In this extension, a color $a$ in an argument means the unit
probability vector $e_a\in\Delta_q$. As in the Ising notation,
we define the extended $\theta_t$ by taking the logarithm of
the right side; its absolute value remains bounded by $J_t^0$.
It is the Boltzmann factor $e^{\theta_t}$ that is affine in each
probability-vector argument separately. We use this
property when prescribing cavity clauses in
Section~\ref{sec:probes}.

\subsection{The M\'ezard--Parisi functional}

The trial spin distributions in the variational formula have
the same hierarchy as in the Ising paper; at each site and state,
we specify a probability vector in place of a magnetization.

We use the rooted tree
$\cA_r=\bigcup_{j=0}^r\N^j$, with root $*$ and terminal vertices
$\N^r$. We write $|\alpha|$ for the depth of a vertex,
$p(\alpha)$ for its path including the root, and
$\alpha\wedge\alpha'$ for the depth of the last common ancestor.
For parameters
\begin{equation}
 0<\zeta_0<\cdots<\zeta_{r-1}<1,
 \label{eq:zetar}
\end{equation}
we denote by $(V_\alpha)_{\alpha\in\N^r}$ the normalized Ruelle
probability cascade \cite{Ruelle}; see \cite[Chapter~2]{SKbook}. 

\begin{definition}[M\'ezard--Parisi trial]\label{def:trial}
A trial consists of a positive integer $r$, $\zz=(\zeta_0,\ldots,\zeta_{r-1})$, and
a measurable function $h_0:[0,1]^r\to\Delta_q$. At every input
site $i$, we sample an independent array of independent uniform variables
$(Z_\gamma^i)_{\gamma\in\cA_r\setminus\{*\}}$, independent also
of the weights $(V_\alpha)_{\alpha\in\N^r}$, and define
\begin{equation}\label{eq:trial}
 m_\alpha^i=h_0\bigl((Z_\gamma^i)_{\gamma\in
                               p(\alpha)\setminus\{*\}}\bigr).
\end{equation}
\end{definition}
The vector $m_\alpha^i$ gives the color probabilities at site $i$
in state $\alpha$. States share the variables on their common
path, whereas different sites use independent arrays. For Ising
spins a local distribution is specified by one magnetization;
here its replacement is the entire vector in $\Delta_q$.

To evaluate a trial, we use the two fields familiar from the
cavity construction: one describes the addition of a site, and
the other corrects the number of clauses among the old sites.
Let $D_t\sim\Poi(k_t\lambda_t)$ and
$\widehat D_t\sim\Poi((k_t-1)\lambda_t)$ be independent counts.
With fresh clauses and a fresh input-site array in every
probability-vector argument, we define
\begin{align}
 A_\alpha^0(a)
 &=\sum_{t\in\mathcal I}\sum_{\nu=1}^{D_t}
       \theta_{t,\nu}(m_\alpha^{t,\nu,1},\ldots,
                     m_\alpha^{t,\nu,k_t-1},a),\label{eq:vertex}\\
 B_\alpha^0
 &=\sum_{t\in\mathcal I}\sum_{\nu=1}^{\widehat D_t}
       \widehat\theta_{t,\nu}(\widehat m_\alpha^{t,\nu,1},\ldots,
                     \widehat m_\alpha^{t,\nu,k_t}).\label{eq:edge}
\end{align}
All counts, clauses, and input-site arrays are mutually independent
and independent of the cascade weights. The counts and clauses
are shared across terminal labels; each input-site array retains the
path sharing in Definition~\ref{def:trial}.
In particular, the hatted clauses and site arrays in the correction field
$B_\alpha^0$ are independent of those in $A_\alpha^0$.
We sum over the color of the added site using the normalized
average $\Av_a f(a)=q^{-1}\sum_{a=1}^q f(a)$. The functional
then compares the two changes in the logarithm of the partition
function.

\begin{definition}[M\'ezard--Parisi functional]\label{def:functional}
For the model \eqref{eq:model}, we define
\begin{equation}\label{eq:Phi}
 \cP_0(r,\zz,h_0)=\log q+
       \E\log\sum_\alpha V_\alpha\Av_a e^{A_\alpha^0(a)}
       -\E\log\sum_\alpha V_\alpha e^{B_\alpha^0}.
\end{equation}
\end{definition}

This is the finite-alphabet, multiple-type version of the functional
in \cite{MPnew}, with the uniform average
over two spins replaced by the uniform average over $q$ colors.
The term $\log q$ restores the counting measure used in the
partition function, so $\cP_0=\log q$ when the interactions
vanish. Assumption~\eqref{eq:primary-moments} makes all the
fields and expectations well defined.

\begin{theorem}[Finite-alphabet M\'ezard--Parisi lower bound]
\label{thm:main}
For the model \eqref{eq:model} under
\eqref{eq:primary-moments},
\begin{equation}\label{eq:formula}
 \liminf_{N\to\infty}F_N
 \ge \inf_{r,\zz,h_0}\cP_0(r,\zz,h_0),
\end{equation}
where the infimum is over the trials in Definition~\ref{def:trial}.
The clause laws may be arbitrary real interactions satisfying
the stated integrability condition.
\end{theorem}

In particular, the lower bound does not require a sign condition
or color-permutation symmetry. A matching upper bound depends on
the interaction. For the antiferromagnetic Potts model, our principal application,
the interpolation in Appendix~\ref{sec:upper} gives this upper bound.

\subsection{The Potts model on random hypergraphs}\label{sec:potts}

For an interaction arity \(K\ge2\), a mean degree \(d\ge0\),
and \(0\le\beta<\infty\), we sample
\(M_N\sim\Poi(dN/K)\) independent ordered \(K\)-tuples of
uniform vertices, with independent entries and allowing repetitions.
The antiferromagnetic Potts Hamiltonian
and pressure are
\begin{equation}\label{eq:potts-model}
 H_N^{K,d,\beta}(\sigma)
 =-\beta\sum_{e=1}^{M_N}
       \one\{\sigma_{i_{e,1}}=\cdots=\sigma_{i_{e,K}}\},
 \qquad
 F_N^{K,d,\beta}
 =\frac1N\E\log\sum_{\sigma\in[q]^N}e^{H_N^{K,d,\beta}(\sigma)}.
\end{equation}
Thus a monochromatic hyperedge is penalized, and the model
has a single primary type with
\[
 k_0=K,\qquad \lambda_0=d/K,\qquad
 \theta_0(a_1,\ldots,a_K)
                   =-\beta\one\{a_1=\cdots=a_K\}.
\]
The graph case is \(K=2\).

We write \(t_\beta=1-e^{-\beta}\) and use a count
\(D\sim\Poi(d)\), with fresh input arrays in every incoming factor.
The extension \eqref{eq:extension} averages the Boltzmann
factor for the monochromatic penalty over the incoming colors.
Summing the resulting site weight over the added color gives
the vertex term below; averaging a single hyperedge gives the edge term.
The cascade identity recalled in Appendix~\ref{sec:upper} separates the independent
clauses in the correction field \eqref{eq:edge}. Consequently \eqref{eq:Phi}
takes the familiar Potts form
\begin{align}
 \mathcal V_\alpha^{K}
 &=\sum_{a=1}^q\prod_{\nu=1}^D
          \left(1-t_\beta\prod_{j=1}^{K-1}
                                    m_{\alpha,a}^{\nu,j}\right),
                                                  \label{eq:potts-vertex}\\
 \mathcal E_\alpha^{K}
 &=1-t_\beta\sum_{a=1}^q\prod_{j=1}^K m_{\alpha,a}^{j},
                                                  \label{eq:potts-edge}\\
 \cP_{K,d,\beta}(r,\zz,h_0)
 &=\E\log\sum_\alpha V_\alpha \mathcal V_\alpha^K
    -\frac{d(K-1)}K\,\E\log\sum_\alpha V_\alpha \mathcal E_\alpha^K.
                                                  \label{eq:potts-Phi}
\end{align}
The incoming clause count has mean \(K\lambda_0=d\);
the correction coefficient is
\((K-1)\lambda_0=d(K-1)/K\).

\begin{theorem}[Hypergraph Potts \MP\ formula]\label{thm:potts}
For every \(q,K\ge2\), \(d\ge0\), and finite \(\beta\ge0\),
\begin{equation}\label{eq:potts-formula}
 \lim_{N\to\infty}F_N^{K,d,\beta}
 =\inf  \cP_{K,d,\beta}(r,\zz,h_0),
\end{equation}
where the infimum is over $r\ge1,$ $\zeta$ in \eqref{eq:zetar}, and
 measurable $h_0:[0,1]^r\to\Delta_q.$
 \end{theorem}
The lower bound is Theorem~\ref{thm:main}. The upper bound follows from the Potts interpolation in
Appendix~\ref{sec:upper}, extending the graph calculation of
\cite{CDGS} to \(K\)-tuples.

The same limit and formula hold for the simple binomial
\(K\)-uniform hypergraph with edge probability
\(d(K-1)!/N^{K-1}\), including \(G(N,d/N)\) when \(K=2\).
The usual coupling changes only \(O(1)\) edges in expectation, for fixed
$K,d,\beta$, with a vanishing effect on the pressure.

\begin{remark}[Symmetry of the overlap matrices]
It remains open whether the infimum in \eqref{eq:potts-formula}
can be restricted to trials whose law is invariant under
permutations of the colors. Let us mention a partial result
in this direction. Given a \MP\ trial in \eqref{eq:trial},
we write
\(
 Q_{\alpha\beta}
   =\E m_\alpha(m_\beta)^{\mathsf T}
\)
for the overlap matrix of two states, where the expectation
integrates the hierarchical coordinates at one site.
This matrix depends only on $\alpha\wedge\beta$. Let
\begin{equation}\label{eq:Psi}
 \Psi(x)=\frac{x}{q}I_q+
       \frac{1-x}{q(q-1)}
          (\one\one^{\mathsf T}-I_q).
\end{equation}
For $x\in[0,1]$, this is the unique matrix of total mass one
and trace $x$ that is invariant under all simultaneous
permutations of its rows and columns. One can prove that the 
infimum in \eqref{eq:potts-formula} can be restricted to finite trials
satisfying the overlap symmetry
\[
 Q_{\alpha\beta}
   =\Psi\bigl(\operatorname{tr}Q_{\alpha\beta}\bigr)
       \qquad\text{for all }\alpha,\beta\in\N^r.
\]
This restriction concerns only the overlap matrices and does
not assert color-permutation invariance of the full trial law.
Since this gives only a partial symmetry reduction, we do
not include the proof here, but briefly indicate the argument.

We combine the symmetry and synchronization argument of
Bates and Sohn \cite{BSPotts} with a self-overlap argument
in the spirit of Chen \cite{ChenSelf}. We first approximate
a small color-symmetric generic Gaussian interaction by weak
diluted clauses at high connectivity. The approximation is
uniform over trials, allowing the general lower bound in
Theorem~\ref{thm:main} to pass to the mixed diluted--Gaussian
model. For the matching upper bound, a further small
antiferromagnetic graph interaction controls the diagonal
interpolation remainder through the $n=1$, $K=2$ term in
\eqref{eq:upper-square}. Differentiability in deterministic
color fields then gives concentration of the self-overlap
at $I_q/q$. Differentiability in the Gaussian coefficients
yields the generalized Ghirlanda--Guerra identities
\cite{GG,PottsMF}; synchronization and color symmetry give
the displayed form of the overlap matrices. We then return
to the diluted approximations and use the reproducibility
construction in this paper to obtain MP trials with
approximate overlap symmetry. After removing the auxiliary
interactions, coarsening the hierarchy and correcting its
covariance increments make the restriction to finite trials
exact.
\end{remark}

\section{The limiting spin law and cavity equations}\label{sec:inputs}

To obtain the lower bound, we use the perturbed limiting spin
law as in the Ising construction. Let us recall its finite-color
form, beginning with the perturbation that supplies the prescribed
clauses later. We may first take finitely many primary clause
types and pass to the countable case at the end; the REM
perturbation itself contains every arity.

\subsection{The diluted REM perturbation}

As in \cite{MPnew}, we use a diluted version of Derrida's random
energy model (REM) \cite{DerridaREM}. Its clauses can approximate
the prescribed cavity interactions in Section~\ref{sec:probes},
while its pressure cost tends to zero with its strength.

For a fixed \(\epsilon>0\) and each arity \(k\ge2\), we sample a Poisson
count \(M_{k,N}\) with mean \(N\) and independent random clauses
\(\theta_{k,\nu}^{\epsilon}:[q]^k\to\mathbb R\), \(1\le\nu\le M_{k,N}\).
Each clause assigns independent Gaussian values to the \(q^k\) possible color configurations:
\[
\theta_{k,\nu}^{\epsilon}(\boldsymbol a)
\sim N(0,\epsilon q^{-4k}),
\qquad \boldsymbol a\in[q]^k.
\]
This independence of the values is what we mean by a REM clause. As in the
original Hamiltonian, we choose the site indices
\(i_{k,\nu,1},\ldots,i_{k,\nu,k}\) independently and uniformly from \([N]\), allowing repetitions, and evaluate the clause at the colors occupying these sites.

We also attach to each site \(i\) a Gaussian vector
$g_i=(g_{i,a})_{a\le q},$ with independent coordinates $g_{i,a}\sim N(0,\epsilon).$ This is a one-site random field: its
contribution to a configuration \(\sigma\) is \(g_{i,\sigma_i}\). All
Poisson counts, clauses, site indices, and Gaussian fields are mutually
independent and independent of the original disorder.

The resulting perturbation is
\[
H_N^{\mathrm{pert},\epsilon}(\sigma)
=
\sum_{k\ge2}\sum_{\nu=1}^{M_{k,N}}
\theta_{k,\nu}^{\epsilon}
\bigl(\sigma_{i_{k,\nu,1}},\ldots,\sigma_{i_{k,\nu,k}}\bigr)
+\sum_{i=1}^{N}g_{i,\sigma_i}.
\]
We add it to \(H_N\), so the perturbed Boltzmann weight is
\(e^{H_N(\sigma)+H_N^{\mathrm{pert},\epsilon}(\sigma)}\).
Let \(F_N^\epsilon\) denote the corresponding pressure.
The variances decay fast enough that the full perturbation is
well defined and has small pressure cost:
\begin{equation}\label{eq:small}
 \sup_N|F_N^\epsilon-F_N|\le C_q\sqrt\epsilon.
\end{equation}
At every arity, the vector of REM clause values has full support
in $\R^{q^k}$. This is what lets us prescribe clauses later.
We keep $\epsilon>0$ throughout the structural argument and
send it to zero only after obtaining the lower bound. The Gaussian
maximum estimates in \cite[Section~7]{GGcavInput} apply here with
$q^k$ clause values and give the required moments for the infinite family.

\subsection{Asymptotic Gibbs measures}

As in \cite{MPnew}, we sample independent replicas and pass to
a subsequential limit of their disorder-averaged spin array.
The Aldous--Hoover representation \cite{Aldous,Hoover}, in the
random-measure formulation of \cite{Austin,Spins}, describes this
limit by a random probability measure $G$ on
\[
 \cH_q=\Bigl\{p\in L^2([0,1];\R^q):p(v)\in\Delta_q
                                      \text{ for almost every }v \Bigr\}.
\]
Conditionally on $G$, we sample independent profiles
$p^1,p^2,\ldots$ with law $G$, and independent uniform
site coordinates $v_1,v_2,\ldots$, independent also of the
profiles and $G$. Given these variables, the spins are
independent and satisfy
\begin{equation}\label{eq:spins}
 \Prob\left(\sigma_i^\ell=a \,\middle|\, G,(p^\ell),(v_i)\right)=p_a^\ell(v_i).
\end{equation}
The independent randomness used to sample these colors is
the terminal, or thermal, randomness in the representation.
Averaging it replaces each indicator
$\one\{\sigma_i^\ell=a\}$ by $p_a^\ell(v_i)$.
Different replicas share the site coordinate $v_i$, so this
average retains their joint color statistics.

In particular, a prescribed color word $\cc=(c_1,\ldots,c_k)$
at one fresh site has conditional probability
\begin{equation}\label{eq:profile-moment}
 \E \left[
      \prod_{\ell=1}^k\one\{\sigma^\ell=c_\ell\}
                  \,\middle|\,G,p^1,\ldots,p^k\right]
       =\int_0^1\prod_{\ell=1}^k p_{c_\ell}^\ell(v)\dd v.
\end{equation}
These are the \emph{colored multioverlaps}. Products of
these moments for independent sites, averaged over the
profiles and $G$, give all finite color-cylinder probabilities.

We write $\langle\cdot\rangle$ for the replica average
under $G^{\otimes\infty}$, and $\E$ for the remaining
randomness. For distinct replicas, their \emph{equality overlap} is defined by
\begin{equation}\label{eq:physical-overlap}
 R_{\ell\ell'}=\int_0^1\sum_{a=1}^q
                    p_a^\ell(v)p_a^{\ell'}(v)\dd v
       =\lim_{M\to\infty}\frac1M\sum_{i=1}^M
                   \one\{\sigma_i^\ell=\sigma_i^{\ell'}\}
                   \quad\text{almost surely}.
\end{equation}
It is the limiting fraction of equal colors, and also the
scalar product of the two profiles in $L^2([0,1];\R^q)$.
The latter interpretation will give the vector metric
\eqref{eq:profile-metric}.

We choose the limiting law to satisfy the spin form of the
Ghirlanda--Guerra identities \cite{GG}:
\begin{equation}\label{eq:GG}
 \E\langle f\psi(R_{1,n+1})\rangle
 =\frac1n\E\langle f\rangle\,\E\langle\psi(R_{12})\rangle
  +\frac1n\sum_{\ell=2}^n\E\langle f\psi(R_{1\ell})\rangle.
\end{equation}
Here $n\ge1$, $f$ is any bounded measurable function of
the spins of the first $n$ replicas, and $\psi$ is any
bounded measurable function of the equality overlap.
Thus $f$ may retain specified colors, even though the new
replica enters through a scalar overlap. The Gaussian and
Poisson regularizations of vanishing pressure cost give
\eqref{eq:GG} together with the cavity equations below,
as in the Ising construction.

To describe the cavity equations for $q$ colors,
we first add one site. Given a state $p$ and the incoming site
coordinates, each incident clause has independent incoming colors
with distributions specified by $p$ at those coordinates.
We average over these colors using the extension
\eqref{eq:extension}, which averages the Boltzmann
factor before taking its logarithm. The argument-position
symmetry imposed in Section~\ref{sec:model} allows us to put
the added color in the last argument. With independent
primary counts $D_t^0\sim\Poi(k_t\lambda_t)$ and REM counts
$D_k\sim\Poi(k)$, we define
\begin{align}
 A(p,a)
 &=g_a+\sum_{t\in\mathcal I}\sum_{\nu\le D_t^0}
   \theta_{t,\nu}(p(v^0_{t,\nu,1}),\ldots,
                         p(v^0_{t,\nu,k_t-1}),a)\notag\\
 &\quad+\sum_{k\ge2}\sum_{\nu\le D_k}
     \theta_{k,\nu}^\epsilon
        (p(v_{k,\nu,1}),\ldots,p(v_{k,\nu,k-1}),a),\label{eq:A}\\
 A(p)&=\log\Bigl(q^{-1}\sum_{a=1}^q e^{A(p,a)}\Bigr),\qquad
 \xi_a(p)=\frac{e^{A(p,a)}}{\sum_b e^{A(p,b)}}.\label{eq:xi}
\end{align}
Summing over the color of the new site multiplies the weight
of state $p$ by $q e^{A(p)}$. The common factor $q$ cancels when
the state weights are normalized, so the state is reweighted
by $e^{A(p)}$. The vector $(\xi_a(p))_{a\le q}\in\Delta_q$
gives the conditional color distribution at that site. In particular,
$\xi_a(p)$ is the expectation of the indicator of color $a$
at the new site, just as $p_a(v)$ is at an original site.
The vector $\xi(p)$ is the counterpart of the Ising cavity
magnetization; its coordinates are probabilities, not sampled
color indicators.

Given $G$, we use independent copies of this cavity construction
at different added sites. For $s$ added sites, the total multiplier
and its normalization are
\[
 \Lambda(p)=\exp\sum_{i=1}^sA_i(p),\qquad
 \mathcal Z=\langle\Lambda\rangle,\qquad
 G^\Lambda(\dd p)=\Lambda(p)G(\dd p)/\mathcal Z.
\]

The cavity equations say that probabilities of prescribed color
assignments are unchanged by this operation, after averaging
the disorder. For $n$ replicas, we specify sets $C_\ell\subset[L]$
of tested sites, where $L\ge s$, and colors $c_{i\ell}$. With this
notation, the identity is
\begin{align}
 &\E\Big\langle\prod_{\ell=1}^n\prod_{i\in C_\ell}
                          p_{c_{i\ell}}^\ell(v_i)\Big\rangle
 \notag\\
 &\quad=\E\int\prod_{\ell=1}^n
      \left(\prod_{i\in C_\ell\cap[s]}\xi_{i,c_{i\ell}}(p^\ell)
      \prod_{i\in C_\ell\setminus[s]}p_{c_{i\ell}}^\ell(v_i)\right)
                         (G^\Lambda)^{\otimes n}(\dd\boldsymbol p).
                                                        \label{eq:cavity}
\end{align}
Thus products of specified coordinates of the original color
vectors become products of the same coordinates of the cavity
outputs. The color labels remain available throughout the identity.

The noncavity coordinates $v_i$ on the right are independent
of the site coordinates used as inputs to the cavity fields.
As in the Ising case, this allows us to include overlap tests:
we may insert the same bounded measurable function $\Psi(R^n)$,
where $R^n=(R_{\ell\ell'})_{\ell<\ell'\le n}$, inside the replica
average on both sides. We use these tests to specify the cluster
genealogy when transferring the cavity identity to barycenters
in Section~\ref{sec:transfer}.

\subsection{The cavity functional and lower bound}

The Aizenman--Sims--Starr cavity
functional \cite{ASS,Spins} in the above model has the same form 
as the M\'ezard--Parisi functional \eqref{eq:Phi}, with the trial
replaced by the asymptotic Gibbs measure $G$ and with the
perturbation included. For the site term we use the log
multiplier $A(p)=\log\Av_a e^{A(p,a)}$ from \eqref{eq:xi}.
The correction field $B(p)$ is the counterpart of
$B_\alpha^0$ in \eqref{eq:edge}: it accounts for the change
in the number of clauses among the old sites when a site
is added. We use independent counts
\(\widehat D_t^0\sim\Poi((k_t-1)\lambda_t)\) and
\(\widehat D_k\sim\Poi(k-1)\). Using fresh clauses and site
coordinates, with all this randomness independent of the random
measure and of that in $A(p,a)$, we define
\begin{align}
 B(p)
 &=\sum_{t\in\mathcal I}\sum_{\nu\le\widehat D_t^0}
    \theta_{t,\nu}(p(v^0_{t,\nu,1}),\ldots,
                                p(v^0_{t,\nu,k_t}))\notag\\
 &\quad+\sum_{k\ge2}\sum_{\nu\le\widehat D_k}
                    \theta_{k,\nu}^\epsilon
                        (p(v_{k,\nu,1}),\ldots,p(v_{k,\nu,k})).
                                                        \label{eq:B}
\end{align}
With these two fields, the cavity functional is
\begin{equation}\label{eq:PG}
 \cP_\epsilon(G)=\log q+
                \E\log\langle e^{A(p)}\rangle
                         -\E\log\langle e^{B(p)}\rangle .
\end{equation}
As in the Ising paper, we can obtain the cavity lower bound
and the two limiting identities for the same law $G$.

\begin{proposition}[Simultaneous finite-alphabet input]\label{prop:input}
For each \(\epsilon>0\), one can add subextensive Gaussian and
Poisson cavity regularizations and select a limiting spin law \(G\)
such that \eqref{eq:GG} and \eqref{eq:cavity} both hold and
\begin{equation}\label{eq:input}
 \liminf_{N\to\infty}F_N^\epsilon\ge\cP_\epsilon(G).
\end{equation}
The regularizations change the pressure by \(o(1)\).
\end{proposition}

\subsection{The joint limiting construction}\label{sec:joint-input}

We can obtain Proposition~\ref{prop:input} from the Ising construction
in \cite[Theorem~1.1]{GGcavInput}, with two adjustments for colors.
Each Gaussian field in the regularization uses products of color
indicators and has covariance equal to the corresponding power
of the equality overlap.
Its integration-by-parts argument still allows every bounded
test of the spin array. In the cavity calculation, we sum over
the $q$ possible colors, obtaining the multiplier and output
probabilities in \eqref{eq:xi}. The moment estimates, finite-volume
comparisons, and joint choice of volumes and regularization
parameters carry over unchanged. Thus both limiting identities
and the lower bound hold for the same law $G$.

We write $\cG_\epsilon$ for the class of random spin laws satisfying
\eqref{eq:GG} and \eqref{eq:cavity} for these fields. From now on
we fix $G\in\cG_\epsilon$. The remaining argument applies to
every law in this class. Our aim is to approximate its spin
distribution by M\'ezard--Parisi trials and thereby bound
$\cP_\epsilon(G)$ from below.

\section{M\'ezard--Parisi reproducibility}\label{sec:geometry}

Proposition~\ref{prop:input} reduces the lower bound to a
comparison of the cavity functional with M\'ezard--Parisi trials.
We now describe the structure of a law $G\in\cG_\epsilon$.
Its states are the profiles $p:[0,1]\to\Delta_q$ introduced above,
with scalar product $R(p,p')$ in $L^2([0,1];\R^q)$. We first
use this scalar product to organize the states into clusters,
then ask what additional information is needed to describe their
joint color probabilities.

\subsection{Overlap geometry and cluster barycenters}

We denote the law of $R_{12}$ under $\E G^{\otimes2}$ by
$\zeta$ and write $q_*=\sup\supp\zeta$. The standard consequences of the
Ghirlanda--Guerra identities apply in the separable Hilbert
space $L^2([0,1];\R^q)$; see \cite{Ultra,SKbook}. They give
\begin{equation}\label{eq:GG-geometry}
 G(\nrm p_{L^2}^2=q_*)=1,\qquad
 R_{12}\ge\min(R_{13},R_{23})\quad\text{almost surely}.
\end{equation}
All these overlaps are nonnegative, since the profiles have
nonnegative coordinates. Here $R_{\ell\ell'}$ is defined for
distinct replicas, as in \eqref{eq:physical-overlap}. The value $q_*$
is the squared norm of a profile; it need not equal
the microscopic self-overlap $1$ of a spin configuration.

The sphere condition in \eqref{eq:GG-geometry} explains why a single
overlap controls the geometry of the whole color vector:
\begin{equation}\label{eq:profile-metric}
 \nrm{p-p'}_{L^2([0,1];\R^q)}^2=2\bigl(q_*-R(p,p')\bigr).
\end{equation}
Thus, two states with overlap close to $q_*$
have close local color distributions in \emph{all colors together},
averaged over the site coordinate. This allows us to approximate
the states by finitely many levels of clusters, as in \cite{MPnew}.
We describe a cut by bins in the overlap range; the top bin
determines how well a cluster can be replaced by its average profile.

\begin{definition}[Finite cut and cluster barycenters]\label{def:cut}
We partition the overlap range into successive bins
$I_j=[t_j,t_{j+1})$ for $j<r$ and $I_r=[t_r,q_*]$,
where $t_0=\inf\supp\zeta$, discarding bins of zero mass
and relabeling the remaining bins in their original order.
An overlap at an internal boundary belongs to the bin on
its right; the boundaries may be atoms.
We measure the two errors of the cut by
\begin{equation}\label{eq:cut-errors}
 \delta=\max_{j<r}\operatorname{diam}(I_j\cap\supp\zeta),
 \qquad \eta=q_*-\inf(I_r\cap\supp\zeta).
\end{equation}
The ultrametric equivalence relations at these boundaries form
a nested cluster tree, whose children we rank by their masses.
For a terminal cluster $C_\alpha$, we write
\[
 V_\alpha=G(C_\alpha),\qquad
 G_\alpha=G(\,\cdot\mid C_\alpha),\qquad
 b_\alpha(v)=\int p(v)G_\alpha(\dd p)\in\Delta_q.
\]
We call $b_\alpha$ its barycenter profile.
\end{definition}

Here $\delta$ measures the unresolved variation in overlaps
below the top level, while $\eta$ measures the size of the
terminal clusters in the metric \eqref{eq:profile-metric}.
We can make both errors tend to zero by allowing the number of
levels to increase. Barycenters approximate the states through
the following bound on the mean squared distance from a state
to the average profile of its cluster:
\begin{equation}\label{eq:bary}
 \int\nrm{p-b_\alpha}_{L^2}^2G_\alpha(\dd p)
       =q_*-\nrm{b_\alpha}_{L^2}^2\le\eta.
\end{equation}
Indeed, the equality is the Hilbert-space variance identity,
and the inequality follows by averaging the scalar product
of two independent states sampled from the same terminal cluster.
The bound holds almost surely for every terminal cluster.
Cauchy--Schwarz then bounds the mean error, averaged over a
uniform site coordinate and a state sampled from that cluster,
by $\sqrt\eta$, also for each specified color coordinate.
This estimate controls the replacement of profiles by
barycenters in the cavity equations.

The geometry between different clusters is retained as well.
If two clusters split at depth $j<r$, the scalar product
of independently sampled states from them lies in
$I_j\cap\supp\zeta$ almost surely. Averaging
shows that the scalar product of their barycenters lies in
its closed convex hull. For a cluster paired with itself,
the corresponding interval is $[q_*-\eta,q_*]$.

\subsection{The joint hierarchical representation}

The cut gives us cluster weights and barycenter profiles. We
now describe their joint law, keeping all color coordinates
together. The cascade parameters are the cumulative masses
of the overlap bins:
\[
 \zeta_j=\zeta(I_0\cup\cdots\cup I_j),\qquad
 0=\zeta_{-1}<\zeta_0<\cdots<\zeta_{r-1}<\zeta_r=1.
\]
The Ghirlanda--Guerra identities identify the cluster weights
with the Ruelle probability cascade having these parameters.
To describe the profiles, we use a stronger consequence of
these identities, hierarchical exchangeability. Here this means
invariance of the vector array under permutations of the cluster
labels that preserve all last-common-ancestor depths.
As in the Ising construction, we use the spin form of the
identities: the tests in \eqref{eq:GG} retain the full color
array while the overlap determines the clusters.

\begin{proposition}[Hierarchical representation at a cut]
\label{prop:heps}
The weights $(V_\alpha)$ have the RPC law above and are
independent of the vector array $(b_\alpha(v_i))_{\alpha,i}$.
This array is hierarchically exchangeable in $\alpha$ and
exchangeable in $i$. It admits a representation in distribution
\begin{equation}\label{eq:he}
 b_\alpha^i=h\bigl((U_\gamma)_{\gamma\in p(\alpha)},
                  (Z_\gamma^i)_{\gamma\in p(\alpha)}\bigr)
       \in\Delta_q,
\end{equation}
where the representing function is deterministic and measurable,
and all the uniform variables are independent of each other
and of the weights. Both families include a root coordinate.
\end{proposition}

We prove the hierarchical exchangeability statement and its
barycenter consequence in Theorem~\ref{thm:heps} and
Corollary~\ref{cor:vector} of Appendix~\ref{app:heps}.
The representation itself is the theorem of Austin and
Panchenko \cite{AP}, applied to a $\Delta_q$-valued array.
This application keeps the color coordinates together and
therefore preserves their joint law.

The two families in \eqref{eq:he} play different roles.
The variables $U_\gamma$ are shared by all sites and can affect
the law of an entire profile. The variables $Z_\gamma^i$ are
independent between sites; at a fixed site, states share the
coordinates on their common path. The M\'ezard--Parisi trial
\eqref{eq:trial} uses only the nonroot coordinates of the second
family. Removing the nonroot $U_\gamma$ is the task of
reproducibility. We retain $U_*$ as a possible global mixture
parameter.

To measure the effect of the shared state variables, we return
to the joint color probabilities at one site. For barycenters,
we use the same products of color probabilities as in
\eqref{eq:profile-moment}.

\begin{definition}[Colored multioverlaps of barycenters]\label{def:colored}
For a tuple of terminal labels
$\boldsymbol\alpha=(\alpha_1,\ldots,\alpha_k)$ and a color
word $\cc=(c_1,\ldots,c_k)\in[q]^k$, we define
\begin{equation}\label{eq:colored-moments}
 M_{\boldsymbol\alpha,\cc}
 =\int_0^1\prod_{\ell=1}^k b_{\alpha_\ell,c_\ell}(v)\dd v
 =\E_Z\prod_{\ell=1}^k b_{\alpha_\ell,c_\ell}^1.
\end{equation}
Here $\E_Z$ integrates all site coordinates of one site,
with the state coordinates held fixed. Each coordinate on a
shared portion of the paths is integrated only once.
Terminal labels may repeat.
\end{definition}

The representation preserves these moments jointly, because
each integral is the almost sure limit of its empirical site
average. Thus the second equality in \eqref{eq:colored-moments}
identifies the original moments with their counterparts in the
represented array, preserving their joint law with the weights.
In particular, summing the two-replica moments over
their common color recovers the scalar product:
$\ip{b_\alpha}{b_{\alpha'}}=\sum_a
M_{(\alpha,\alpha'),(a,a)}$. We may choose a representation
separately at each cut; all comparisons concern the spin
statistics it represents at that cut.

\subsection{Shapes and approximate reproducibility}

We can now formulate the property needed to remove the
nonroot state variables. Once the genealogy of a finite set
of replicas and the global variable $U_*$ are fixed, a colored
multioverlap should have little remaining randomness. We first
record that genealogy in a form that keeps track of the color
assigned to each replica.

The \emph{labeled shape} $T$ of $k$ terminal labels is their
matrix of last-common-ancestor depths. An entry equals $r$
when the corresponding terminal labels coincide. The row and
column labels identify the replica occurrences, so the color
word assigns a color to each of them even when terminal labels
repeat. We retain these labels in the moment identification
of Section~\ref{sec:graph}. To record the corresponding overlaps,
we choose strictly increasing representatives
\[
 q_j\in\overline{\operatorname{conv}}(I_j\cap\supp\zeta)
       \quad(j<r),\qquad q_r=q_*.
\]
We associate with $T$ its quantized overlap matrix
\begin{equation}\label{eq:quantized}
 q(T)_{\ell\ell'}=
 \begin{cases}
 q_{\alpha_\ell\wedge\alpha_{\ell'}},&\ell\ne\ell',\\
 q_*,&\ell=\ell'.
 \end{cases}
\end{equation}
Almost surely, it approximates both the scalar products of the
cluster barycenters and the overlaps of sampled states in these
clusters entrywise to within $\delta\vee\eta$.
The diagonal comparison uses the squared profile norms,
not microscopic spin self-overlaps.

Let $\E_U$ integrate the nonroot state coordinates in the
union of the replica paths, retaining $U_*$. By
\eqref{eq:he}, the conditional mean
$\E_U[M_{\boldsymbol\alpha,\cc}\mid U_*]$ depends only on
$U_*$, the labeled shape, and the color word. We measure the
remaining dependence by
\begin{equation}\label{eq:D}
 D_{k,\cc}=\E\sum_{\boldsymbol\alpha}
       V_{\alpha_1}\cdots V_{\alpha_k}
 \left(M_{\boldsymbol\alpha,\cc}
       -\E_U[M_{\boldsymbol\alpha,\cc}\mid U_*]\right)^2.
\end{equation}
Thus $D_{k,\cc}$ is the average conditional variance of a
colored multioverlap given its shape and the global variable,
when terminal labels are sampled independently conditional on
the weights. Independence of the weights and the represented
array makes this an average over the RPC shape law.
For the finite-alphabet analogue of the Ising defect, we take
$D_k=\max_{\cc\in[q]^k}D_{k,\cc}$. This maximum matters
when a child branch inherits only part of a color word: we
will need the lower-order estimate for whichever colors occur
there. As in the Ising case, the defect is zero for a
M\'ezard--Parisi field with a fixed global parameter.

\begin{theorem}[Approximate reproducibility]\label{thm:repro}
Let $\epsilon>0$ and $G\in\cG_\epsilon$. Along every sequence
of finite cuts with $\delta,\eta\to0$,
\[
 \max_{\cc\in[q]^k}D_{k,\cc}\longrightarrow0
       \qquad\text{for each fixed }k.
\]
The conclusion holds for every choice of hierarchical
representation at each cut, with no restriction on the
numbers of overlap levels.
\end{theorem}

The absence of a depth restriction lets us refine a general
overlap distribution, rather than assume finitely many overlap
values. The cavity equations supply the additional control
needed for this conclusion; the geometry and representation
above use only the Ghirlanda--Guerra identities.

\section{Sketch of proof}\label{sec:proof}

We follow the induction in \cite{MPnew}, replacing its nonnegative
Ising variables $S=(1+b)/2$ by the color probabilities $b_c$.
We recall the quantities needed to describe this replacement and
give the adjustments to the argument. Throughout this section,
$G\in\cG_\epsilon$, $\epsilon>0$, and the number of replicas are
fixed, while $\delta,\eta\to0$. All estimates are uniform in the
retained depth and the number of overlap levels.

\subsection{The induction for colored moments}\label{sec:norm}

We begin with the hierarchical representation \eqref{eq:he},
retaining a common path through depth $p<r$. We write
$U_{\le p},Z_{\le p}$ for its state and site coordinates,
including the root. For a
labeled $k$-replica shape $T$, let $U_T,Z_T$ denote the future
coordinates in the union of its paths. A shared coordinate occurs
only once. The analogue of the future multioverlap in the Ising case is the future colored multioverlap
\begin{equation}\label{eq:future}
 F_{T,\cc}^p=\E_{Z_T}\prod_{\ell=1}^k b_{\alpha_\ell,c_\ell},
                    \qquad \cc\in[q]^k.
\end{equation}
All state variables remain fixed in this expectation. Thus
$F_{T,\cc}^p$ is the conditional probability of the prescribed
color word after averaging the future site coordinates.

As in the Ising case, to prove approximate reproducibility, we can use the covariance norm of \cite[Definition~6.2]{MPnew}.
Writing $U$ for the future state variables and $Z$ for the
retained site coordinates, let
\begin{equation}\label{eq:norm}
 \Pi_p f=f-\E_Uf,\qquad
 \nrm f_p^4=\E_{U_{\le p}}\iint
       \bigl(\E_U f(U,Z)f(U,Z')\bigr)^2\dd Z\dd Z'.
\end{equation}
The state prefix is suppressed from the arguments of $f$.
For a centered function, the inner expectation is the covariance,
over the future state variables, between the evaluations at two
independent site prefixes. 

Let $\mu_k^p$ be the shape law of the auxiliary RPC
with parameters $\zeta_p,\ldots,\zeta_{r-1}$, as in the Ising
induction. The stronger induction quantity is now
\begin{equation}\label{eq:kappa}
 \kappa_k=\max_{\cc\in[q]^k}\sup_{p<r}
       \sum_T\mu_k^p(T)\nrm{\Pi_pF_{T,\cc}^p}_p^4.
\end{equation}
Here we also take the maximum over color words $\cc\in[q]^k$,
so the induction controls every choice of colors.
As in \cite[Lemma~6.4]{MPnew},
$D_{k,\cc}\le\sqrt{\kappa_k}$, so it suffices to prove
$\kappa_k\to0$.

For the base case $k=1$, we work on one path and write
$F_a(U,Z)=\E_{Z_{>p}}b_a(U_{\le p},U,Z,Z_{>p})$. Centering in
the future state variables gives
\begin{equation}
 f_a=F_a-\E_UF_a,\qquad
 C_{ab}(Z,Z')=\E_Uf_a(U,Z)f_b(U,Z').
 \label{eq:covab}
\end{equation}
If we apply the geometric argument of \cite[Lemma~6.3]{MPnew} to $F(U,\cdot)$ as a vector in $L^2(\dd Z;\R^q)$, the resulting estimates become
\begin{align}
 \E\int\sum_a|f_a(U,Z)|^2\dd Z&\le\delta,
                                                  \label{eq:geom-L2}\\
 \E_{U_{\le p}}\sum_{a,b}\iint|C_{ab}(Z,Z')|^2\dd Z\dd Z'
                              &\le\delta^2/4.     \label{eq:geom-HS}
\end{align}
Thus the scalar
overlap controls every color coordinate and every pair of colors
separately, without cancellation between colors. The diagonal
terms $a=b$ give $\kappa_1\le\delta^2/4$. We will also use the
control of the full covariance matrix in the cavity equations below.

For the induction step, we assume $\kappa_m\to0$ for $m<k$.
As in the Ising case, we obtain the common-stem moment $H_{T,\cc}^p$
from $F_{T,\cc}^p$ by averaging the future state variables in every
child branch after the first split $a=a(T)$, retaining the common stem. This is the same
common-stem moment as in \cite[Lemma~6.6]{MPnew}. If all
terminal labels coincide, there are no such branches to average.
The branch-replacement estimate in \cite{MPnew} now gives
\begin{equation}\label{eq:branch-bound}
 \max_{\cc}\sup_{p<r}\sum_T\mu_k^p(T)
       \nrm{F_{T,\cc}^p-H_{T,\cc}^p}_p^4
       \le k^4B_k\max_{m<k}\kappa_m,
\end{equation}
where $B_k$ is the number of set partitions of $[k]$.
Its proof is unchanged: a child containing the replica
indices $D$ carries the restricted color word $(c_\ell)_{\ell\in D}$,
and the maximum over colors in the definition \eqref{eq:kappa} supplies the corresponding
lower-order estimate.

\subsection{Color covariances under the cavity change of density}
\label{sec:transfer}

The cavity equations \eqref{eq:cavity} pass to the cluster
barycenters and their RPC representation just as in
\cite[Section~7]{MPnew}. Using the log multiplier $A_i(b_\alpha)$
from \eqref{eq:xi}, we write the terminal log multiplier for $s$
added sites on one canonical path as
$X_r=\sum_{i=1}^s A_i(b_\alpha)$. If $\omega_j$ collects its
state coordinate and all incoming site coordinates at level $j$,
the backward recursion and path density \cite{RPC,SKbook} are
\begin{equation}\label{eq:rho}
 X_j=\zeta_j^{-1}\log\E_{\omega_{j+1}}e^{\zeta_jX_{j+1}},
 \qquad
 \rho_{p:r}=\prod_{j=p}^{r-1}e^{\zeta_j(X_{j+1}-X_j)}.
\end{equation}
This change of
density depends on the full incoming color vectors: the sum
over output colors has already been taken in each $A_i$.
It is the same density for every color observable.

For two cavity sites, let $\Cov^{12}_{p\to r}$ denote covariance
over the future state and input-site coordinates under their
joint density above, with the prefix and cavity disorder fixed.
The geometric estimate \eqref{eq:geom-HS} transfers to
\begin{equation}\label{eq:cov-color}
 \sup_{p<r}\E
   \sum_{a,b}
    \bigl|\Cov^{12}_{p\to r}(\xi_{1,a},\xi_{2,b})\bigr|^2
                    \longrightarrow0.
\end{equation}
Here $\xi_{i,a}$ is the probability of color $a$ at cavity
site $i$, as in \eqref{eq:xi}. The expectation in
\eqref{eq:cov-color} includes the cavity disorder and the
common prefix under the RPC change of density.

Let us explain why the scalar equality overlap suffices
for this transfer. We choose four paths $\alpha_1,\ldots,\alpha_4$ that agree through depth $p$ and enter distinct children immediately afterwards.
Let
\[
R_{\ell m}^{a}
:=\int b_{\alpha_\ell,a}(v)b_{\alpha_m,a}(v)\,dv,
\qquad
R_{\ell m}:=\sum_a R_{\ell m}^{a}.
\]
As in the Ising case, expanding the square of the covariance kernel $C_{ab}(Z,Z')$ in \eqref{eq:covab} using these four conditionally independent continuations, and then averaging the two site coordinates, gives
\[
 \E_{U_{\le p}}\iint |C_{ab}(Z,Z')|^2\dd Z\dd Z'
 =
 \E\left[
 R_{12}^{a}R_{12}^{b}
 -2R_{12}^{a}R_{13}^{b}
 +R_{13}^{a}R_{24}^{b}
 \right].
\]
We do not know the individual colored overlaps $R_{\ell m}^{a}$, but summing over $a,b$ gives
\[
\begin{aligned}
 \E_{U_{\le p}}\sum_{a,b}
       \iint |C_{ab}(Z,Z')|^2\dd Z\dd Z'
 &=
 \E\left[R_{12}^2-2R_{12}R_{13}+R_{13}R_{24}\right]\\
 &=
 \E\left[(R_{12}-R_{13})(R_{12}-R_{24})\right]
 \le\delta^2.
\end{aligned}
\]
The second equality uses exchangeability of the four
branches. The inequality follows because all the scalar
products appearing here lie in the same interval of width
at most $\delta$.

To transfer this estimate, we apply the cavity identity
\eqref{eq:cavity} to the underlying products before averaging
the two site coordinates. For example, the first term
comes from
\(
 b_{\alpha_1,a}^{1}b_{\alpha_1,b}^{2}
 b_{\alpha_2,a}^{1}b_{\alpha_2,b}^{2},
\)
and the identity replaces each factor by the corresponding
cavity-output coordinate, with its color unchanged.
As in \cite[Propositions~7.1--7.2]{MPnew}, the same
three-term combination on the cavity side is the square
of the conditional covariance under the joint path density.
The errors from replacing states by barycenters tend to
zero uniformly in $p$. Since every color pair uses the
same density, summing over $a,b$ proves \eqref{eq:cov-color}.
The nonnegativity of these squared covariances then gives
control of each color pair separately.

\subsection{Prescribed clauses and the mixed cavity tests}
\label{sec:probes}

The full support of the REM clauses and the Gaussian one-site
fields allows us to prescribe cavity interactions, by the same
localization of the nonnegative covariance estimate as in
\cite[Proposition~8.1]{MPnew}. Let us describe the clauses we need.
For a positive separately affine function
$\mathcal A:(\Delta_q)^m\to(0,\infty)$, we prescribe
\begin{equation}\label{eq:prescribed-clause}
 e^{\vartheta(a_1,\ldots,a_m,c)}
 =\begin{cases}
   \mathcal A(e_{a_1},\ldots,e_{a_m}),&c=1,\\
   1/(q-1),&c\ne1,
  \end{cases}
\end{equation}
and choose the one-site field to be $\log z$ at color $1$
and zero at the other colors, where $z>0$.
The extension \eqref{eq:extension} gives weight $z\mathcal A$
to color $1$ and total weight $1$ to the remaining colors.
With input spin $\chi(c)=2\one\{c=1\}-1$, we therefore have
\begin{equation}\label{eq:probe}
 e^{A}=\frac{1+z\mathcal A}{q},\qquad
 \xi^\chi=2\xi_1-1=\frac{z\mathcal A-1}{z\mathcal A+1}.
\end{equation}
This is the Ising multiplier and response (or cavity output), up to a constant
in the multiplier. One output color carries $\mathcal A$,
which may depend on every color of every input vector.
To relate these scalar outputs to the color probabilities, recall
that $\xi_i^\chi=2\xi_{i,1}-1$ at cavity site $i$. Under the joint
cavity change of density, we therefore have
\[
 \Cov^{12}_{p\to r}(\xi_1^\chi,\xi_2^\chi)
 =4\,\Cov^{12}_{p\to r}(\xi_{1,1},\xi_{2,1}).
\]
The covariance on the right is the color-$1$, color-$1$ term
controlled by \eqref{eq:cov-color}. We thus obtain the same
convergence to zero in mean square for the conditional covariance
of the scalar outputs, uniformly in the retained depth $p$.

From here through the mixed cavity tests, the proofs in
\cite[Sections~8--9]{MPnew} apply without change. In particular, if
$\rho^1,\rho^2$ are the densities \eqref{eq:rho} for the
terminal fields $\log(1+z_i\mathcal A_i)$ and $\rho^{12}$
is the density for their sum, the covariance control gives
$\sup_{p<r}\|\rho_{p:r}^{12}-\rho_{p:r}^1\rho_{p:r}^2\|_2\to0$
under the original product law.
These superscripts index cavity sites, with separate input
arrays and a common state path; each density uses all colors.
The resulting one-site estimate, the marked-replica expansion,
and the extraction of a mixed coefficient need no further
adjustment. They use the same scalar response
\eqref{eq:probe}, while the separately affine input functions
now have probability-vector arguments.

We recall the resulting test in the notation needed below.
Let $P_1,\ldots,P_k:(\Delta_q)^m\to\R$ be bounded separately
affine functions, and write
$P_v^{\alpha_\ell}=P_v(b_{\alpha_\ell}^1,\ldots,b_{\alpha_\ell}^m)$.
Let $\Sp_T$ average the future state variables outside replica
$1$'s path, retaining that path and all site prefixes.
The mixed cavity test is
\begin{equation}\label{eq:mixed}
 \sup_{p<r}\left\|
   \sum_T\mu_k^p(T)\Sp_T\Pi_p\E_{Z_T}
       \sum_{\pi\in S_k}\prod_{v=1}^kP_v^{\alpha_{\pi(v)}}
                 \right\|_p\longrightarrow0.
\end{equation}
This is \cite[Proposition~9.3]{MPnew}.
The permutation sum assigns each $P_v$ to a different
replica in the mixed coefficient. All functions use the same
input sites; different sites have independent site arrays,
and replicas share coordinates on their common paths.

\subsection{Identifying each colored moment}\label{sec:graph}

Given \eqref{eq:mixed}, we now use the multigraph tests of
\cite[Section~10]{MPnew}. We fix a color word $\cc\in[q]^k$
and a loopless multigraph $\Gamma$ on the labeled vertices
$[k]$. Writing $q_{v,w}(T)$ for the entries of the overlap
matrix $q(T)$, we associate with $\Gamma$ the monomial
\[
 P_\Gamma(q(T))
   =\prod_{\{v,w\}\in E(\Gamma)}q_{v,w}(T),
\]
where edges are counted with their multiplicities.

In the test function associated with vertex $v$, we replace
the distinguished Ising factor $S_0$ by the color coordinate
$x^0_{c_v}$ at input site $0$. We retain one additional input
site for each edge and sum over the common color in the two
factors associated with that edge. Averaging the edge-site
coordinates gives scalar products of the corresponding
barycenters, which we approximate by the entries of $q(T)$
as in the Ising proof.

In the permutation sum in \eqref{eq:mixed}, the function
associated with vertex $v$ is assigned to replica $\pi(v)$.
We write $T^\pi$ for the shape of the reordered tuple
$(\alpha_{\pi(1)},\ldots,\alpha_{\pi(k)})$ and set
$c^\pi_\ell=c_{\pi^{-1}(\ell)}$. Thus
\[
 P_\Gamma(q(T^\pi))
   =\prod_{\{v,w\}\in E(\Gamma)}
                      q_{\pi(v),\pi(w)}(T),
\]
while the factor at input site $0$ gives
$F_{T,\cc^\pi}^p$. The calculation in
\cite[Section~10]{MPnew} therefore gives
\begin{equation}\label{eq:before-branch}
 \sup_{p<r}\left\|
  \sum_T\mu_k^p(T)\sum_{\pi\in S_k}
     P_\Gamma(q(T^\pi))\,\Sp_T\Pi_pF_{T,\cc^\pi}^p
                    \right\|_p\longrightarrow0.
\end{equation}
The color word now depends on $\pi$, so the future
multioverlap remains inside the permutation sum. We handle
this dependence by first replacing it with its common-stem
moment.

Recall that $\Pi_pH_{T,\cc^\pi}^p$ depends only on future
state variables on the common stem. Since this stem belongs
to every replica path, $\Sp_T$ leaves the centered moment
unchanged:
\(
 \Sp_T\Pi_pH_{T,\cc^\pi}^p=\Pi_pH_{T,\cc^\pi}^p.
\)
The contraction properties of $\Sp_T$ and $\Pi_p$ consequently
give
\[
 \bigl\|\Sp_T\Pi_pF_{T,\cc^\pi}^p
                -\Pi_pH_{T,\cc^\pi}^p\bigr\|_p
 \le \bigl\|F_{T,\cc^\pi}^p-H_{T,\cc^\pi}^p\bigr\|_p.
\]
Since $P_\Gamma$ is bounded on the overlap matrices,
\eqref{eq:branch-bound} allows
us to make this replacement in \eqref{eq:before-branch} uniformly over color words, with error tending to zero uniformly in $p$ by induction.

We can now relabel the replicas by $\beta_v=\alpha_{\pi(v)}$. This simply rewrites the same
product of color coordinates,
\(
 \prod_{\ell=1}^k b_{\alpha_\ell,c^\pi_\ell}
   =\prod_{v=1}^k b_{\beta_v,c_v},
\)
in terms of replicas with shape $T^\pi.$
After averaging the branch coordinates, the common stem
is unchanged, and we obtain the equality
\(
 H_{T,\cc^\pi}^p=H_{T^\pi,\cc}^p
\)
as functions of the same canonical stem coordinates.
In particular, the choice of the retained replica path
no longer affects the moment.

Finally, we interchange the sums over shapes and permutations.
For each fixed $\pi$, the map $T\mapsto T^\pi$ is a bijection
of the labeled shapes, and replica exchangeability gives
$\mu_k^p(T)=\mu_k^p(T^\pi)$. Using $T'=T^\pi$ as the summation
variable, we obtain
\begin{align}\label{eq:relabel}
 &\sum_T\mu_k^p(T)\sum_{\pi\in S_k}
       P_\Gamma(q(T^\pi))\Pi_pH_{T^\pi,\cc}^p
=\sum_{\pi\in S_k}\sum_{T'}\mu_k^p(T')
       P_\Gamma(q(T'))\Pi_pH_{T',\cc}^p
       \notag\\
 &\qquad=k!\sum_T\mu_k^p(T)
       P_\Gamma(q(T))\Pi_pH_{T,\cc}^p.
\end{align}
Thus each permutation gives the same sum after changing
the shape variable. We obtain the original graph monomial
with the color word $\cc$ fixed, while the common-stem
moment retains its dependence on the labeled shape and
the assigned colors.

Polynomial approximation now gives
\begin{equation}\label{eq:continuous-test}
 \sup_{p<r}\left\|
       \sum_T\mu_k^p(T)\psi(q(T))\Pi_pH_{T,\cc}^p
                     \right\|_p\longrightarrow0
\end{equation}
for every continuous function $\psi$ of the labeled overlap
matrix. We retain the replica labels here because different
replicas may carry different colors. This ordering of branch
replacement and relabeling is the adjustment to the Ising
identification argument.

\subsection{Completion of the argument}\label{sec:compact}
\label{sec:bridge}

The compactness argument of \cite[Section~11]{MPnew} now
applies to \eqref{eq:continuous-test}, keeping the overlap
cells labeled. The only change in its moment recursion is
at the leaves, where $S$ becomes $b_{c_\ell}\in[0,1]$.
Its martingale estimate for moving branching depths already
allows different nonnegative child functions, by polarization,
so both that estimate and its independence of the number of
levels are unchanged. The remaining measure estimates use
only the scalar overlap identities. The closing argument in
\cite[Section~12]{MPnew} therefore gives $\kappa_k\to0$,
using \eqref{eq:branch-bound} and taking the maximum over
the finitely many color words. This proves
Theorem~\ref{thm:repro}.

Finally, colored multioverlaps determine the finite spin
distributions, as in \eqref{eq:profile-moment}. The resampling
and variational passage of \cite[Section~13]{MPnew} therefore
give
\[
 \cP_\epsilon(G)\ge
       \inf_{r,\zz,h_0}\cP_\epsilon(r,\zz,h_0),
\]
where the trial functional includes the perturbation.
Its cost is $O_q(\sqrt\epsilon)$ uniformly over trials,
as for the pressure in \eqref{eq:small}. Proposition~\ref{prop:input}
and $\epsilon\downarrow0$ prove Theorem~\ref{thm:main}.
The same truncation of primary clause types, justified by
\eqref{eq:primary-moments}, gives the countable case.

\appendix
\section{The Potts interpolation upper bound}\label{sec:upper}

We interpolate between the Potts model and a fixed trial,
using the same three types of terms as in the graph calculation
of \cite{CDGS}.
We first recall the cascade identity that evaluates
the two ends of the interpolation; see \cite{RPC,SKbook}.
If $X_\alpha$ is
a bounded hierarchical field generated by independent uniforms
along paths, with these uniforms and any global disorder independent
of the cascade weights, we define the backward recursion
\[
 X_j=\frac1{\zeta_j}\log\E_j e^{\zeta_jX_{j+1}},
                      \qquad j=r-1,\ldots,0,
\]
where \(X_r\) is the terminal value and \(\E_j\) integrates
the next coordinate given the prefix and the global disorder. Then
\begin{equation}\label{eq:rpc-log}
 \E\log\sum_\alpha V_\alpha e^{X_\alpha}=\E X_0.
\end{equation}
The outer expectation includes any global disorder held fixed
in the backward recursion. If several fields use independent
coordinate arrays conditional on this disorder, their exponential integrals
factor at every level, and the resulting $X_0$ for their sum is the sum of
their individual $X_0$'s. We may therefore separate the
independent clauses in a correction field by first conditioning
on the counts and clauses. For a Poisson
number of independent identically distributed clauses, independent
of that count, the expected cascade logarithm is
the mean count times the expected cascade logarithm for one
clause. Truncation extends this calculation to the integrable
fields in \eqref{eq:Phi}.

We keep the trial fixed and let $u\in[0,1]$ be the interpolation parameter.
As the parameter increases, we replace one-site trial factors
by model hyperedges and add the corresponding trial correction
factors. The endpoint comparison gives the upper bound
if the interpolated pressure is nonincreasing.
We define a partition function over
$(\sigma,\alpha)\in[q]^N\times\N^r$, with weight $V_\alpha$,
using the following three independent families of factors:
\begin{enumerate}[label=(\roman*),itemsep=2pt]
\item \(\Poi(udN/K)\) model hyperedge factors
  \(
    1-t_\beta\one\{\sigma_{i_1}=\cdots=\sigma_{i_K}\};
  \)
  with independent uniform vertex indices, allowing repetitions;
\item at each vertex \(i\), \(\Poi((1-u)d)\) one-site factors
  \(
    1-t_\beta\prod_{j=1}^{K-1}m_{\alpha,\sigma_i}^{j},
  \)
  with \(K-1\) fresh independent input arrays for each factor;
\item \(\Poi(ud(K-1)N/K)\) trial correction factors
  \(
    1-t_\beta\sum_{a=1}^q\prod_{j=1}^K m_{\alpha,a}^{j},
  \)
  with \(K\) fresh independent arrays for each factor.
\end{enumerate}
All counts, vertex indices, and input arrays are independent
of the cascade weights. Each input array uses the fixed trial
function and is shared across terminal labels according to their paths.
Let \(\varphi_N(u)\) be the expected logarithm of this
partition function divided by \(N\). Independence and
\eqref{eq:rpc-log} give
\begin{align}
 \varphi_N(0)&=\E\log\sum_\alpha V_\alpha \mathcal V_\alpha^K,
                                                   \label{eq:upper-left}\\
 \varphi_N(1)&=F_N^{K,d,\beta}
       +\frac{d(K-1)}K\E\log\sum_\alpha V_\alpha \mathcal E_\alpha^K.
                                                   \label{eq:upper-right}
\end{align}

To determine the sign of the derivative, we compare the color
patterns produced by the actual spins and the trial \eqref{eq:trial} at a fresh site. We sample $n$ replicas $(\sigma^\ell,\alpha_\ell)$ independently from the interpolating Gibbs measure, conditional on its disorder. For a color word
$\boldsymbol a\in[q]^n$, we write
\[
 R_{\boldsymbol a}
 =\frac1N\sum_{i=1}^N
                  \prod_{\ell=1}^n\one\{\sigma_i^\ell=a_\ell\},
 \qquad
 M_{\boldsymbol a}(\boldsymbol\alpha)
 =\E \prod_{\ell=1}^n
                                      m_{\alpha_\ell,a_\ell}.
\]
The first quantity is the empirical frequency of the color
word $\boldsymbol a$ among the $N$ sites. The second is its
probability at a fresh trial site, after integrating the whole
hierarchical input array with the replica labels held fixed.
This fresh array is independent of the interpolating disorder;
replicas use the same coordinates on common portions of their paths.
Both are probability vectors when indexed by
$\boldsymbol a\in[q]^n$.

Poisson differentiation and the expansion
\(\log(1-t_\beta x)=-\sum_{n\ge1}t_\beta^n x^n/n\) give
\begin{equation}\label{eq:upper-square}
 \varphi_N'(u)
 =-\frac dK\sum_{n\ge1}\frac{t_\beta^n}{n}
     \E\Big\langle\sum_{\boldsymbol a\in[q]^n}
       \bigl(R_{\boldsymbol a}^{K}
          -K R_{\boldsymbol a}M_{\boldsymbol a}^{K-1}
          +(K-1)M_{\boldsymbol a}^{K}\bigr)
                                             \Big\rangle_u\le0.
\end{equation}
Here brackets denote the conditional replica average in the
interpolating system, and the outer expectation averages its disorder.
The three terms come from the three families of factors:
the model hyperedge contributes
\(\sum_{\boldsymbol a}R_{\boldsymbol a}^{K}\);
the one-site factor contributes
\(\sum_{\boldsymbol a}R_{\boldsymbol a}M_{\boldsymbol a}^{K-1}\);
and the correction contributes
\(\sum_{\boldsymbol a}M_{\boldsymbol a}^{K}\).
The rates at which these three families change give the
coefficients $1,-K,K-1$. Their combination is nonnegative
because $x\mapsto x^K$ is convex on $[0,\infty)$:
\[
 x^K-Kxy^{K-1}+(K-1)y^K
     =x^K-y^K-Ky^{K-1}(x-y)\ge0.
\]
For $K=2$, this remainder is the familiar square $(x-y)^2$;
the same convexity argument works for every $K\ge2$ because
the quantities being compared are nonnegative. Each of the
three sums over color words lies in $[0,1]$, which also
justifies the series expansion at every finite temperature,
where $t_\beta<1$.

Combining \eqref{eq:upper-left}--\eqref{eq:upper-square}
gives \(F_N^{K,d,\beta}\le\cP_{K,d,\beta}(r,\zz,h_0)\)
for every trial. Together with Theorem~\ref{thm:main}, this
proves Theorem~\ref{thm:potts}. The cases \(d=0\) or \(\beta=0\)
give \(\log q\) directly. We handle finite mixtures of these
antiferromagnetic hyperedge interactions by adding their
nonnegative convexity remainders.

\section{A short proof of hierarchical exchangeability}
\label{app:heps}

We include this appendix because the Dirichlet calculation below
gives a new short proof of the hierarchical exchangeability
theorem in \cite{HEPS}. We present it in the finite-alphabet form
used in Proposition~\ref{prop:heps}, keeping the full color array
together. Throughout the appendix, we assume only the spin
Ghirlanda--Guerra identities \cite{GG} in the form \eqref{eq:GG}.

\subsection{The cut and its sampled tree}

We use the cut and clusters of Definition~\ref{def:cut}, with
$\zeta=\Law(R_{12})$ and $q_*=\sup\supp\zeta$. As in
\eqref{eq:GG-geometry}, the profiles lie on the sphere
$\nrm p^2=q_*$ and their scalar products are ultrametric.
Writing the cut thresholds as $u_1<\cdots<u_r$, we have
\begin{equation}\label{eq:zetas}
 \zeta_{j-1}=\zeta([0,u_j)),\qquad
 0=\zeta_{-1}<\zeta_0<\cdots<\zeta_{r-1}<\zeta_r=1.
\end{equation}
The relation $R(p,p')\ge u_j$ defines the depth-$j$
clusters. We enumerate the children of each cluster by decreasing
mass, writing $C_\alpha$ for the cluster, $V_\alpha=G(C_\alpha)$
for its mass, and $G_\alpha$ for its normalized restriction. The weights have the
corresponding RPC law. In particular the sibling masses
are positive and distinct almost surely. The point-mass
overlap law gives a one-state measure and can be handled
directly, so we consider nontrivial cuts.

The finite-alphabet extension of \cite{HEPS} follows by
carrying the full spin array as a mark in the
Ghirlanda--Guerra prediction rule. The following
Dirichlet calculation isolates the independence statement
that we need before passing to ranked clusters.

We sample \(n\ge1\) profiles \(p^1,\ldots,p^n\) from \(G\)
and write \(S^n=(\sigma_i^\ell)_{i\ge1,\,\ell\le n}\) for their
spin array, sampled as in \eqref{eq:spins}.
To describe how these profiles fall into the clusters, we form
a finite tree \(T\). At each depth \(j=0,\ldots,r\), its vertices
are the depth-\(j\) clusters containing at least one sampled
profile, with edges joining each cluster to its visited children.
A vertex may have only one child; for example, a single sampled
profile gives a path of length \(r\).
Each leaf at depth \(r\) carries the labels of the replicas
belonging to that cluster. We fix a labeled tree shape \(T\) and let \(O_T\) be the event
that the sampled replicas \(p^1,\ldots,p^n\) form this shape,
with the specified replica labels at each leaf.
We consider only shapes for which \(\Prob(O_T)>0\).

For \(t\in T\), let \(C_t\) be its cluster; \(C_*\) is the whole
space. For a nonroot node \(t\), we use one original replica
descending from it as a representative. Its overlap with a new profile determines
membership in \(C_t\), independently of which such representative
was chosen. We partition the whole space into the residual regions
\begin{equation}\label{eq:regions}
 B_t=
 \begin{cases}
 C_t\setminus\displaystyle\bigcup_{u\text{ occupied child of }t}C_u,
                                                    & |t|<r,\\
 C_t,                                               & |t|=r.
 \end{cases}
 \qquad \delta_t=G(B_t).
\end{equation}
The finite vector \(\delta=(\delta_t)_{t\in T}\) has nonnegative
coordinates summing to one. It records both the masses of the
occupied terminal clusters and the as yet unvisited parts of their
ancestors. The occupied cluster masses are recovered by
\[
 V_t=\sum_{\substack{u\in T\\u\text{ descendant of }t}}\delta_u,
\]
where a node is counted as its own descendant.

Let \(b_t\) be the number of occupied children of an internal node,
and \(n_t\ge1\) the number of original replicas in a leaf. We define
\begin{equation}\label{eq:dir-params}
 a_t=
 \begin{cases}
 b_t\zeta_{|t|}-\zeta_{|t|-1},& |t|<r,\\
 n_t-\zeta_{r-1},& |t|=r.
 \end{cases}
\end{equation}
The root case is \(a_*=b_*\zeta_0\).
All parameters are strictly positive: every occupied internal node
has at least one child, and \(\zeta_j>\zeta_{j-1}\).
The terms associated with edges cancel in the sum, giving
\begin{equation}\label{eq:dir-sum}
 \sum_{t\in T}a_t=\sum_{|t|=r}n_t=n.
\end{equation}

\begin{theorem}[Dirichlet law of residual cluster masses]\label{thm:dirichlet}
Under the spin Ghirlanda--Guerra identities \eqref{eq:GG}, conditional on \(O_T\), the residual
mass vector \(\delta\) is independent of \(S^n\), and
\begin{equation}\label{eq:dir-law}
 \Law(\delta\mid S^n,O_T)=\Dir((a_t)_{t\in T}).
\end{equation}
Equivalently, if \(m_t\) are nonnegative integers, \(M=\sum_t m_t\),
and \((a)_m=a(a+1)\cdots(a+m-1)\), with \((a)_0=1\), then
\begin{equation}\label{eq:dir-moment}
 \E\left\langle f(S^n)\one_{O_T}\prod_{t\in T}\delta_t^{m_t}
       \right\rangle
 =\E\langle f(S^n)\one_{O_T}\rangle
            \frac{\prod_{t\in T}(a_t)_{m_t}}{(n)_M}.
\end{equation}
In particular \(S^n\) and the occupied cluster masses
\((V_t)_{t\in T}\) are independent conditional on \(O_T\).
\end{theorem}

\begin{proof}
We keep the regions in \eqref{eq:regions} fixed from the original
\(n\) replicas as we draw further replicas. After \(j\) such draws,
we write \(K_t(j)\) for the number of new replicas lying in \(B_t\).
For a nonroot node \(t\) at depth \(h\), let \(N_t(j)\) be the number
of all \(n+j\) replicas in \(C_t\).

We apply the Ghirlanda--Guerra identities \eqref{eq:GG}, choosing
the old representative of \(C_t\) as the distinguished replica
and using \(\psi(x)=\one\{x\ge u_h\}\). Since the tests include
every bounded function of \(S^{n+j}\), we obtain the following
conditional prediction on \(O_T\):
\begin{equation}\label{eq:cluster-predict}
 \Prob(p^{n+j+1}\in C_t\mid S^{n+j}, O_T)
 =\frac{1-\zeta_{h-1}+N_t(j)-1}{n+j}
 =\frac{N_t(j)-\zeta_{h-1}}{n+j}.
\end{equation}
Every membership event appearing here is measurable from the
overlaps with old representatives, hence from \(S^{n+j}\).

For an internal node, we subtract \eqref{eq:cluster-predict} for all
its occupied children from the corresponding formula for the node. There were no original replicas in its
residual region, so
\(N_t(j)-\sum_{u\text{ child}}N_u(j)=K_t(j)\).
At the root, we subtract the child probabilities from \(1\).
For a leaf, \(N_t(j)=n_t+K_t(j)\).
All three cases give the single formula
\begin{equation}\label{eq:urn}
 \Prob(p^{n+j+1}\in B_t\mid S^{n+j},O_T)
       =\frac{a_t+K_t(j)}{n+j}.
\end{equation}
Thus future replicas enter these fixed regions according to a
P\'olya urn with initial parameters \(a_t\). In particular, this
prediction contains no additional dependence on the original
color-spin marks.

For a sequence \(t_1,\ldots,t_M\) in which region \(t\) appears
\(m_t\) times, we multiply the successive conditional probabilities
\eqref{eq:urn}. Keeping an arbitrary \(f(S^n)\one_{O_T}\) in the
expectation gives
\begin{align*}
 &\E\left\langle f(S^n)\one_{O_T}
       \prod_{j=1}^M\one\{p^{n+j}\in B_{t_j}\}\right\rangle
 =\E\langle f(S^n)\one_{O_T}\rangle
                       \frac{\prod_t(a_t)_{m_t}}{(n)_M}.
\end{align*}
On the other hand, conditional on \(G\), the original profiles,
and the site variables and spin randomness of the first \(n\)
replicas, future profiles are independent with law \(G\).
Integrating them out changes the product of indicators into
\(\prod_t\delta_t^{m_t}\). This proves \eqref{eq:dir-moment}.

The displayed moments are those of the Dirichlet distribution;
for example, this follows by integrating its density proportional
to \(\prod_t x_t^{a_t-1}\) on the simplex.
Polynomials determine probability measures on a compact simplex.
Applying this observation first with \(f=1\), and then with every
bounded \(f\), proves both the law and the conditional independence.
\end{proof}

The passage from \eqref{eq:urn} to the Dirichlet law is the
classical P\'olya-urn calculation; see \cite{BlackwellMacQueen}.
For one overlap level, the conditional law of the visited
cluster masses and the remaining mass, given the sampled
partition, is a special case of \cite[Corollary~20]{PitmanUrn},
with parameters $\alpha=\zeta_0$ and $\theta=0$.
Here we derive the prediction rule directly from the spin
Ghirlanda--Guerra identities \eqref{eq:GG}. By keeping an
arbitrary bounded function of $S^n$ throughout the calculation,
we also obtain the independence of the residual masses from
the sampled spin array, conditional on the occupied tree.

\subsection{From sampled trees to ranked clusters}\label{sec:ranked}

The rest of the argument closely follows the final section, Section~5, of \cite{HEPS}, and the proof below is a somewhat compressed version of it.

Theorem~\ref{thm:dirichlet} concerns clusters visited by
finitely many replicas. To obtain Proposition~\ref{prop:heps},
we pass to clusters ranked by their true
masses. We first note that unused replicas can be removed
without changing the conditional law of the retained marks.

\begin{lemma}[Deleting unused replicas]\label{lem:delete}
Let \(T\) be a finite labeled overlap-tree shape and \(I\) any
subset of its replica labels. Let \(T|_I\) be the restricted
shape, retaining the original depth levels. Then the conditional
law of the spin marks of the retained replicas, given \(O_T\),
is their conditional law given \(O_{T|_I}\).
\end{lemma}

\begin{proof}
By replica exchangeability, we can delete the last replica and
then repeat. Suppose there are \(N\) replicas before this last
one is added. Its prescribed relation to them has one of two forms.
It either joins an occupied terminal cluster with \(s\) replicas,
or first enters a new child at depth \(a+1\) of an occupied
depth-\(a\) cluster.

In the first case \eqref{eq:cluster-predict} gives probability
\((s-\zeta_{r-1})/N\). In the second case, let the parent contain
\(s\) replicas and have \(b\) occupied children with counts
\(s_1,\ldots,s_b\), where \(\sum_i s_i=s\).
The required event is to join the parent and avoid \emph{all}
its occupied children. Its conditional probability is
\begin{equation}\label{eq:delete-prob}
 \frac{s-\zeta_{a-1}}N
       -\sum_{i=1}^b\frac{s_i-\zeta_a}N
       =\frac{b\zeta_a-\zeta_{a-1}}N.
\end{equation}
For the root, the first fraction is \(1\), since
\(s=N\) and \(\zeta_{-1}=0\).
Below that new child there are no old replicas, so no further
conditions on overlaps with old replicas are imposed.

In either case, the extension probability is a positive constant
determined by the old shape, independent of all its spin marks.
Multiplying by any bounded function of those marks and dividing
by the shape probability proves the assertion for one deletion.
Iterating this deletion proves the lemma. Note that the subtraction of all occupied
children in \eqref{eq:delete-prob} accounts for the entire event
that the replica enters a new child.
\end{proof}

We now sample independently in every terminal cluster,
conditional on \(G\):
\[
 p^{\alpha,1},p^{\alpha,2},\ldots\ \sim G_\alpha,
       \qquad \alpha\in\N^r.
\]
We use common independent site coordinates \(v_i\) and sample
conditionally independent colors \(S_i^{\alpha,\ell}\) with probabilities
\(p_a^{\alpha,\ell}(v_i)\).
A permutation of \(\N^r\) is \emph{hierarchy preserving} if it
preserves the depth of the last common ancestor of every pair.
Equivalently, it independently permutes children at the vertices
of the rooted tree.

\begin{theorem}[Hierarchical exchangeability of the color array]\label{thm:heps}
Assume \eqref{eq:GG}. For every finite cut, the entire marked array
\(
 (S_i^{\alpha,\ell})_{\alpha\in\N^r,\ i,\ell\ge1}
\)
is independent of the cluster-mass array
\((V_\gamma)_{\gamma\in\cA_r}\). Its law is invariant under
hierarchy-preserving permutations of \(\alpha\), under site
permutations, and under separate permutations of the sample
index \(\ell\) in each terminal cluster.
\end{theorem}

\begin{proof}
For fixed \(d,L,M\), we consider the first \(d\) children at every
vertex through depth \(r\), the first \(L\) samples in each of
their terminal clusters, and the first \(M\) sites.
We approximate this array using \(n\) ordinary independent samples
from \(G\). Among the children visited by these samples, we rank
by their \emph{true} masses and select the first \(d\), recursively,
then retain the first \(L\) visits to each selected terminal cluster.
We call the procedure successful when all required clusters and
visits exist.

For a fixed labeled sample shape \(T\), the ranking, success event,
and selected replica labels are functions of the occupied masses
and \(T\). Theorem~\ref{thm:dirichlet} makes those masses independent
of the full sampled spin array conditional on \(T\).
For any particular successful selection of labels,
Lemma~\ref{lem:delete} says that their conditional spin law is
determined only by their restricted shape: a \(d\)-ary tree of
depth \(r\), with \(L\) replicas in each leaf. We denote this fixed
law on the retained \(M\)-site marks by \(K_{d,L,M}\).
It does not depend on the larger shape or the masses used to
perform the selection.

Consequently, for every bounded function \(F\) of the selected
marks and every bounded function \(H\) of the selected masses,
\begin{equation}\label{eq:finite-factor}
 \E[FH\one_{\mathrm{success}}]
       =K_{d,L,M}(F)\,\E[H\one_{\mathrm{success}}].
\end{equation}
We obtain this equality by summing over the finitely many labeled
shapes \(T\) and selected replica-label lists. On each part, the
selection event and \(H\) depend only on the mass vector, whereas
the conditional expectation of \(F\) is the same constant.

Almost surely each of the finitely many true leading clusters
has positive mass. Eventually it is visited, and every selected
terminal cluster receives at least \(L\) visits. Once the true
first \(d\) children have been visited, their order among the
visited children is correct. Thus the procedure succeeds and
stabilizes to the intended ranked clusters and their first
\(L\) visits as \(n\to\infty\).
Conditional on \(G\), the subsequence of visits to each cluster
is an independent sample from its normalized restriction; these
subsequences are independent across clusters. This is the
within-cluster sampling used in the statement.
Bounded convergence in \eqref{eq:finite-factor} proves independence
for the specified finite arrays.

The balanced restricted shape defining \(K_{d,L,M}\) is invariant
under tree automorphisms of the finite \(d\)-ary tree. Replica
exchangeability therefore makes \(K_{d,L,M}\) invariant under
those automorphisms. Any two finite lists of leaf labels with
the same ancestor relations fit in such a tree after increasing
\(d\). This proves hierarchical exchangeability.
Site exchangeability and the within-cluster replica permutations
follow from the sampling construction.
Finally, we let \(d,L,M\) grow. Since finite cylinders determine
the laws of both infinite arrays, this proves the theorem.
\end{proof}

\subsection{The joint barycenter representation}
\label{sec:representation}

We return to the barycenter $b_\alpha(v)=\int p(v)G_\alpha(\dd p)$
from Definition~\ref{def:cut}. The preceding theorem gives its
joint representation with the cluster weights.

\begin{corollary}[Joint vector representation]\label{cor:vector}
The array $(b_\alpha(v_i))_{\alpha,i}$ is independent of
the cluster masses, hierarchically exchangeable in
$\alpha$, and exchangeable in $i$. It has a representation
\begin{equation}\label{eq:Hexch}
 b_\alpha^i\ \stackrel{\mathrm d}=
 h\bigl((U_\gamma)_{\gamma\in p(\alpha)},
        (Z_\gamma^i)_{\gamma\in p(\alpha)}\bigr)\in\Delta_q,
\end{equation}
with independent uniform variables, independent also of
the masses. The representation preserves jointly all
colored multioverlaps and scalar products of barycenters.
\end{corollary}

\begin{proof}
For every cluster, site, and color, the conditional law
of large numbers gives simultaneously
\begin{equation}\label{eq:freq}
 b_{\alpha,a}(v_i)=\lim_{L\to\infty}\frac1L
       \sum_{\ell=1}^L\one\{S_i^{\alpha,\ell}=a\}.
\end{equation}
Thus the vector array inherits the independence and
exchangeability in Theorem~\ref{thm:heps}. Applying the
hierarchical representation theorem of \cite{AP} to this
$\Delta_q$-valued array, we obtain \eqref{eq:Hexch} with the two
families of variables in \eqref{eq:he}: the $U$'s are shared
across sites, whereas each
site has its own $Z$ array. Both families include the root.

Finally, for each tuple of labels and colors,
\[
 \int_0^1\prod_{\ell=1}^k b_{\alpha_\ell,c_\ell}(v)\dd v
 =\lim_{M\to\infty}\frac1M\sum_{i=1}^M
                     \prod_{\ell=1}^k b_{\alpha_\ell,c_\ell}(v_i)
 \quad\text{almost surely}.
\]
There are countably many such tuples. The equality of
the vector-array laws therefore carries all these
quantities jointly to the representation. Summing the
two-replica moments over their common color gives the
scalar products as well.
\end{proof}

This proves Proposition~\ref{prop:heps}: the scalar overlap
organizes the clusters, while the full-spin tests preserve their
joint color information.

\end{document}